\documentclass[11pt]{article}
\usepackage[labelsep=period]{caption}
\usepackage{bm}
\usepackage{framed}
\usepackage{fullpage}
\usepackage{epsfig}
\usepackage{graphics}
\usepackage{latexsym}
\usepackage{amsmath}
\usepackage{amsfonts}
\usepackage{amssymb}
\usepackage{mathrsfs}
\usepackage{pifont}
\usepackage{yhmath}
\usepackage{undertilde}
\usepackage{bbm}
\usepackage{dsfont}
\usepackage{eufrak}
\usepackage{wasysym}
\usepackage{underscore}
\usepackage{epstopdf}
\usepackage{fontenc}
\usepackage{amsthm}
\newtheorem{theorem}{Theorem}[section]

\newtheorem{conjecture}{Conjecture}[section]

\def\qed{\hfill \rule{4pt}{7pt}}
\newcommand{\de}{\backslash}
\DeclareMathAlphabet{\mathpzc}{OT1}{pzc}{m}{it}
\title{\bf On Gupta's Co-density Conjecture}
\author{\vspace{2mm} 
Yan Cao$^{a}$ \quad Guantao Chen$^{a}$ \quad Guoli Ding$^{b}$ \quad Guangming Jing$^{a}$\thanks{E-mail: gjing1@gsu.edu.} \quad 
Wenan Zang$^{c}$\thanks{G. Chen was supported in part by NSF grant DMS-1855716 and  NSFC grant 11871239.
G. Ding was supported in part by NSF grant DMS-1500699. W. Zang was supported in part by the Research Grants 
Council of Hong Kong.}\\
$\stackrel{a}{}$ Department of Mathematics and Statistics, Georgia State University\\
Atlanta, GA 30303, USA \smallskip\\
$\stackrel{b}{}$ Mathematics Department, Louisiana State  University\\ Baton Rouge, LA 70803, USA \smallskip\\
$\stackrel{c}{}$ Department of Mathematics, The University of Hong Kong\\Hong Kong, China }
\begin{document}
\date{}
\maketitle

\begin{abstract}
Let $G=(V,E)$ be a multigraph. The {\em cover index} $\xi(G)$ of $G$ is the greatest integer $k$ for which there is
a coloring of $E$ with $k$ colors such that each vertex of $G$ is incident with at least one edge of each color.
Let $\delta(G)$ be the minimum degree of $G$ and let $\Phi(G)$ be the {\em co-density} of $G$, defined by  
\[\Phi(G)=\min \Big\{\frac{2|E^+(U)|}{|U|+1}:\,\, U \subseteq V, \,\, |U|\ge 3 \hskip 2mm 
{\rm and \hskip 2mm odd} \Big\},\]
where $E^+(U)$ is the set of all edges of $G$ with at least one end in $U$. It is easy to see that $\xi(G) \le 
\min\{\delta(G), \lfloor \Phi(G) \rfloor\}$. In 1978 Gupta proposed the following co-density conjecture:
Every multigraph $G$ satisfies $\xi(G)\ge \min\{\delta(G)-1, \, \lfloor \Phi(G) \rfloor\}$, which is
the dual version of the Goldberg-Seymour conjecture on edge-colorings of multigraphs. In this note we prove that 
$\xi(G)\ge \min\{\delta(G)-1, \, \lfloor \Phi(G) \rfloor\}$ if $\Phi(G)$ is not integral and 
$\xi(G)\ge \min\{\delta(G)-2, \, \lfloor \Phi(G) \rfloor-1\}$ otherwise. We also show that this co-density conjecture 
implies another conjecture concerning cover index made by Gupta in 1967.
\end{abstract}

\openup 0.5\jot

\newpage

\section{Introduction}

In this note we consider multigraphs, which may have parallel edges but contain no loops. Let $G=(V,E)$ be
a multigraph. The {\em chromatic index} $\chi'(G)$ of $G$ is the least integer $k$ for which there is
a coloring of $E$ with $k$ colors such that each vertex of $G$ is incident with at most one edge of 
each color. Let $\Delta(G)$ be the maximum degree of $G$ and let $\Gamma(G)$ be the {\em density} of 
$G$, defined by
\[\Gamma(G)=\max \Big\{\frac{2|E(U)|}{|U|-1}:\,\, U \subseteq V, \,\, |U|\ge 3 \hskip 2mm
{\rm and \hskip 2mm odd}\Big\},\]
where $E(U)$ is the set of all edges of $G$ with both ends in $U$.  Clearly, $\chi'(G) \ge \max\{\Delta(G),
\, \Gamma(G) \}$; this lower bound, as shown by Seymour \cite{Se} using Edmonds' matching polytope theorem 
\cite{e65}, is precisely the {\em fractional chromatic index} of $G$, which is the optimal value of the 
{\em fractional edge-coloring problem}:
\begin{center}
\begin{tabular}{ll}
\hbox{Minimize} \ \ \ & $ {\bm 1}^T {\bm x}$  \\
\hbox{\hskip 0.2mm subject to} & $ A{\bm x} = {\bm 1}$  \\
& \hskip 3mm ${\bm x} \ge {\bm 0},$
\end{tabular}
\end{center}
where $A$ is the edge$-$matching incidence matrix of $G$. In the 1970s Goldberg \cite{G73} and Seymour \cite{Se}
independently made the following conjecture.
\begin{conjecture}\label{GS}
Every multigraph $G$ satisfies $\chi'(G)\le \max\{\Delta(G)+1, \, \lceil \Gamma(G) \rceil\}$.
\end{conjecture}

Over the past four decades this conjecture has been a subject of extensive research, and has stimulated an
important body of work, with contributions from many researchers; see McDonald \cite{Mc} for a survey  
on this conjecture and Stiebitz {\em et al.} \cite{SSTF} for a comprehensive account of edge-colorings.   
Recently, three of the authors, Chen, Jing, and Zang, have announced a complete
proof of Conjecture \ref{GS} \cite{CJZ}.

The present note is devoted to the study of the dual version of the classical edge-coloring problem (ECP), 
which asks for a coloring of the edges of $G$ using the maximum number of colors in such a way that at each vertex 
all colors occur. It is easy to see that each color class induces an edge cover of $G$. (Recall that an {\em edge cover} 
is a subset $F$ of $E$ such that each vertex of $G$ is incident to at least one edge in $F$.) So this problem is actually 
the {\em edge cover packing problem} (ECPP). Let $\xi(G)$ denote the optimal value of ECPP, which we call the 
{\em cover index} of $G$. As it is $NP$-hard \cite{H} in general to determine the chromatic index $\chi'(G)$
of a simple cubic graph $G$, determining the cover index $\xi(G)$ is also $NP$-hard.

Let $\delta(G)$ be the minimum degree of $G$, let $E^+(U)$ be the set of all edges of $G$ with at least one end in 
$U$ for each $U \subseteq V$, and let $\Phi(G)$ be the {\em co-density} of $G$, defined by
\[\Phi(G)=\min \Big\{\frac{2|E^+(U)|}{|U|+1}:\,\, U \subseteq V, \,\, |U|\ge 3 \hskip 2mm 
{\rm and \hskip 2mm odd} \Big\}.\]
Obviously, $\xi(G) \le \delta(G)$.  Since each edge cover contains at least $(|U|+1)/2$ edges in $E^+(U)$ for any
$U\subseteq V$ with $|U| \ge 3$ and odd, $\Phi(G)$ provides another upper bound for $\xi(G)$. So
$\xi(G) \le \min\{\delta(G), \Phi(G)\}$. Based on a polyhedral description of edge covers (see Theorem 27.3 in 
Schrijver \cite{Sc}), Zhao, Chen, and Sang \cite{ZCS} observed that the parameter $\min\{\delta(G), \, \Phi(G)\}$ 
is exactly the {\em fractional cover index} of $G$, the optimal value of the {\em fractional edge cover packing 
problem} (FECPP):
\begin{center}
\begin{tabular}{ll}
\hbox{Maximize} \ \ \ & $ {\bm 1}^T {\bm x}$  \\
\hbox{\hskip 0.3mm subject to} & $ B{\bm x} = {\bm 1}$  \\
& \hskip 3mm ${\bm x} \ge {\bm 0},$  
\end{tabular}
\end{center}
where $B$ is the edge$-$edge cover incidence matrix of $G$.  They \cite{ZCS} also devised a
combinatorial polynomial-time algorithm for finding the co-density $\Phi(G)$ of any multigraph $G$. 

In 1978 Gupta \cite{G2} proposed the following co-density conjecture, which is the counterpart of 
Conjecture \ref{GS} on ECPP. 
\vskip 2mm
\begin{conjecture}\label{CC1}
Every multigraph $G$ satisfies $\xi(G)\ge \min\{\delta(G)-1, \, \lfloor \Phi(G) \rfloor\}$.
\end{conjecture}

The reader is referred to Stiebitz {\em et al.} \cite{SSTF} for more information about this conjecture.
Its validity would imply that, first, there are only two possible values for the cover 
index $\xi(G)$ of a multigraph $G$: $\min\{\delta(G)-1, \, \lfloor \Phi(G) \rfloor\}$ and $\min\{\delta(G), \, 
\lfloor \Phi(G) \rfloor\}$; second, any multigraph has a cover index within one of its fractional cover 
index, so FECPP also has a fascinating integer rounding property (see  Schrijver \cite{Sc86,Sc}); third, 
even if $P \ne NP$, the $NP$-hardness of ECPP does not preclude the possibility of designing an efficient 
algorithm for finding at least $\min\{\delta(G)-1, \, \lfloor \Phi(G) \rfloor\}$ disjoint edge covers 
in any multigraph $G$.  

To our knowledge, the bound $\xi(G) \ge \min\{\lfloor \frac{7 \delta(G)+1}{8} \rfloor, \lfloor 
\Phi(G) \rfloor\}$ established by Gupta \cite{G2} in 1978 remains to be the best approximate version 
of Conjecture \ref{CC1}.  

As is well known,  the inequality $\chi'(G)\le \Delta(G)+ \mu(G)$ holds for any multigraph $G$, where $\mu(G)$ 
is the maximum multiplicity of an edge in $G$. This result has been successfully dualized by Gupta \cite{G} to 
packing edge covers: $\xi(G)\ge \delta(G)-\mu(G)$. It is worthwhile pointing out that this dual version follows 
from Conjecture \ref{CC1} as a corollary, because $\Phi(G)\ge \delta(G)-\mu(G)$. To see this, let $U$ be a subset of 
$V$ with $|U|\ge 3$ and odd, let $F(U)$ be the set of all edges of $G$ with precisely one end in $U$, and 
let $G[U]$ be the subgraph of $G$ induced by $U$. Since each vertex in $U$ is adjacent to at most $(|U|-1) \mu(G)$
edges in $G[U]$ and at most $|F(U)|$ edges outside $G[U]$, we have $\delta(G)\le (|U|-1) \mu(G)+ |F(U)|$, which 
implies that $\delta(G)|U|+|F(U)| \ge  (\delta(G)-\mu(G))(|U|+1)$. As $2|E^+(U)| =2|E(U)|+2|F(U)|\ge \delta(G)|U|
+|F(U)|$, we obtain $2|E^+(U)| \ge (\delta(G)-\mu(G))(|U|+1)$ and hence $\Phi(G) \ge \delta(G)-\mu(G)$, as desired.

Gupta \cite{G} demonstrated that the lower bound $\delta(G)-\mu(G)$ for $\xi(G)$ is sharp when $\mu(G)\ge 1$ and $\delta(G)=
2p\mu(G)-q$, where $p$ and $q$ are two integers satisfying $q\ge 0$ and $p > \mu(G)+ \lfloor (q-1)/2 \rfloor$. This 
led Gupta \cite{G} to suggest the following conjecture, which aims to give a complete characterization of all values of 
$\delta(G)$ and $\mu(G)$ for which no multigraph $G$ with $\xi(G) = \delta(G)-\mu(G)$ exists. 
\vskip 2mm
\begin{conjecture}\label{CC2}
Let $G$ be a multigraph such that $\delta(G)$ cannot be expressed in the form $2p\mu(G)-q$, for any two integers $p$ and $q$
satisfying $q\ge 0$ and $p > \mu(G)+\lfloor (q-1)/2 \rfloor$. Then $\xi(G)\ge \delta(G)-\mu(G)+1$.
\end{conjecture}

As edge covers are more difficult to manipulate than matchings, it is no surprise that a direct proof of 
conjecture \ref{CC1} would be more complicated and sophisticated than that of Conjecture \ref{GS} (see 
\cite{CJZ}, which is under review). One purpose of this note is to establish a slightly weaker version of conjecture 
\ref{CC1} by using Conjecture \ref{GS}.     

\begin{theorem} (Assuming Conjecture \ref{GS}) \label{cover} 
Let $G$ be a multigraph. Then $\xi(G)\ge \min\{\delta(G)-1, \, \lfloor \Phi(G) \rfloor\}$ if $\Phi(G)$ is not 
integral and $\xi(G)\ge \min \{\delta(G)-2, \, \lfloor \Phi(G) \rfloor-1 \}$ otherwise. 
\end{theorem}

{\bf Remark}. Suppose $\Phi(G)<\delta(G)$. By this theorem, we obtain $\xi(G)=\lfloor\Phi(G)\rfloor$ if $\Phi(G)$ 
is not integral and $\Phi(G)-1\le \xi(G) \le \Phi(G)$ otherwise, because $\xi(G)\le \min\{\delta(G), \, \lfloor 
\Phi(G) \rfloor\}$.  

\vskip 2mm
In this note we also show that Conjecture \ref{CC2} is contained in Conjecture \ref{CC1} as a special case. 

\begin{theorem} \label{bound} 
Conjecture \ref{CC1} implies Conjecture \ref{CC2}.
\end{theorem}

Throughout this note we shall repeatedly use the following terminology and notations. Let $G=(V,E)$ be a multigraph. A
subset $U$ of $V$ is called an {\em odd set} if $|U|$ is odd and $|U|\ge 3$. For each $v \in V$, let 
$d_G(v)$ be the degree of $v$ in $G$. For each $U \subseteq V$, let $E_G(U)$ be the set of all edges of $G$ with 
both ends in $U$, let $E_G^+(U)$ be the set of all edges of $G$ with at least one end in $U$, and let $F_G(U)$ 
be the set of all edges of $G$ with exactly one end in $U$.  For any two subsets $X$ and $Y$ of $V$, let 
$E_G(X,Y)$ be the set of all edges of $G$ with one end in $X$ and the other end in $Y$. We write $E_G(x,y)$ 
for $E_G(X,Y)$ if $X=\{x\}$ and $Y=\{y\}$. We shall drop the subscript $G$ if there is no danger of confusion.

The proofs of the above two theorems will take up the entire remainder of this note.  

\section{Approximate Version} 

We present a proof of Theorem \ref{cover} in this section. Let $G=(V,E)$ be a multigraph and let $Z\subseteq V$. A 
set $C\subseteq E$ is called a {\it $Z$-cover} if every vertex of $Z$ is incident with at least one edge of $C$. Note 
that if $Z=V$, then $Z$-covers are precisely edge covers of $G$. Let $e\in E(x,y)$ and let $G'$ be obtained from
$G$ by adding a new vertex $x'$ and making $e$ incident with $x'$ instead of $x$ (yet still incident with 
$y$); we say that $G'$ arises from $G$ by {\em splitting off} $e$ from $x$.  To prove the theorem, we shall 
actually establish the following variant.

\begin{theorem} \label{cover100} 
Let $G=(V,E)$ be a multigraph, let $Z\subseteq V$, let $k$ be a positive integer, and let $\epsilon$ be $0$ or $1$.
If $d(z)\ge k+1$ for all $z\in Z$ and $|E^+(U)|\ge \frac{|U|+1}{2}k + \epsilon$ for all odd sets $U\subseteq Z$, 
then $G$ contains $k-1+\epsilon$ disjoint $Z$-covers. 
\end{theorem} 
 
{\bf Proof}. Splitting off edges from vertices outside $Z$ if necessary, we may assume that all vertices outside $Z$ 
have degree one. Suppose for a contradiction that Theorem \ref{cover100} is false. We reserve the triple $(G,Z, k)$ for 
a counterexample with the minimum $\sum_{z\in Z} d(z)$. For convenience, we call an odd set $U \subseteq Z$ {\em optimal} 
if $|E^+(U)|= \frac{|U|+1}{2}k+ \epsilon$.

By hypothesis, $d(z)\ge k+1$ for all $z\in Z$, which can be strengthened as follows.  

\vskip 2mm
{\bf Claim}.  $d(z)=k+1$ for all $z\in Z$. 

\vskip 2mm

Otherwise, $d(z)\ge k+2$ for some $z\in Z$.  If $z$ is contained in no optimal odd set $U \subseteq Z$, 
letting $H$ be obtained from $G$ by splitting off an edge from $z$, then $(H,Z,k)$ would be a smaller counterexample 
than $(G,Z, k)$, a contradiction. Hence 

(1) there exists an optimal odd set $U_1 \subseteq Z$ containing $z$; subject to this, we assume that $|U_1|$ is minimum.

Since $(|U_1|+1)k+ 2 \epsilon =2|E^+(U_1)|=2|E(U_1)|+2|F(U_1)| \ge (k+1)|U_1|+|F(U_1)|$, we have $|F(U_1)|\le k-|U_1|
+ 2 \epsilon \le k < d(z)$.
So $z$ is adjacent to some vertex $y \in U_1$. Let $H$ be arising from $G$ by splitting off one edge 
$e \in E(y,z)$ from $z$. We propose to show that 

(2) $(H,Z,k)$ is a smaller counterexample than $(G,Z, k)$.
 
Assume the contrary. Then $|E^+_H(U_2)|< \frac{|U_2|+1}{2}k+ \epsilon$ for some odd set $U_2 \subseteq Z$ by the hypothesis
of this theorem. Thus

(3) $z\in U_2$, $y \notin U_2$, and $|E^+(U_2)|=\frac{|U_2|+1}{2}k+ \epsilon$. 

Let $T_1=U_1\de U_2$ and  $T_2= U_2\de U_1$. By (3), we have $y\in U_1\de U_2$, so $T_1 \ne \emptyset$. By
the minimality assumption on $|U_1|$ (see (1)), $U_2$ is not a proper subset of $U_1$, which implies
$T_2 \ne \emptyset$. Since $z\in U_1\cap U_2$, we obtain $|U_1\cap U_2|\ge 1$. Let us consider two cases, 
according to the parity of $|U_1\cap U_2|$.

{\bf Case 1.} $|U_1\cap U_2|$ is odd.

It is a routine matter to check that 

(4) $|E^{+}(U_1\cup U_2)|+|E^{+}(U_1\cap U_2)|=|E^{+}(U_1)|+|E^{+}(U_2)|-|E(T_1, T_2)|$.

In this case, $U_1 \cup U_2$ is an odd set. So $|E^+(U_1\cup U_2)|\ge \frac{|U_1 \cup U_2|+1}{2}k + \epsilon$ by the hypothesis
of this theorem. 

(5) $|E^+(U_1\cap U_2)|\ge \frac{|U_1 \cap U_2|+1}{2}k+\epsilon+1$.

To justify this, note that if $|U_1\cap U_2|=1$, then $|E^+(U_1\cap U_2)|=d(z) \ge k+2$. So (5) holds. 
If $|U_1\cap U_2|\ge 3$, then  $U_1\cap U_2$ is not an optimal odd set by the minimality assumption on $|U_1|$ 
(see (1)). Thus (5) is also true. 

From (4) and (5) we deduce that $\frac{|U_1 \cup U_2|+1}{2}k + \epsilon \le |E^{+}(U_1\cup U_2)| \le 
|E^{+}(U_1)|+|E^{+}(U_2)| - |E^{+}(U_1\cap U_2)| \le \frac{|U_1|+1}{2}k + \epsilon + \frac{|U_2|+1}{2}k + \epsilon
- \frac{|U_1 \cap U_2|+1}{2}k-\epsilon-1 = \frac{|U_1 \cup U_2|+1}{2}k + \epsilon -1$, a contradiction.

{\bf Case 2.} $|U_1\cap U_2|$ is even.

It is easy to see that $|E^{+}(U_1)|+|E^{+}(U_2)|=|E^{+}(T_1)|+|E^{+}(T_2)|+2|E(U_1\cap U_2)|+|E(U_1\cap 
U_2,T_1\cup T_2)|+ 2|E(U_1\cap U_2,\overline{U_1\cup U_2})|$, where $\overline{U_1\cup U_2}=V- (U_1\cup U_2)$.
Thus

(6) $|E^{+}(U_1)|+|E^{+}(U_2)|\ge |E^{+}(T_1)|+|E^{+}(T_2)|+2|E(U_1\cap U_2)|+|F(U_1\cap U_2)|$.

In this case, $|T_i|$ is odd, so $|E^{+}(T_i)| \ge \frac{|T_i|+1}{2}k+ \epsilon$ for $i=1,2$ by the hypothesis of 
this theorem. It follows from (3) and (6) that $\frac{|U_1|+1}{2}k + \epsilon + \frac{|U_2|+1}{2}k + \epsilon
\ge \frac{|T_1|+1}{2}k + \epsilon + \frac{|T_2|+1}{2}k + \epsilon +  2|E(U_1\cap U_2)|+|F(U_1\cap U_2)| \ge \frac{|T_1|+1}{2}k 
+ \epsilon + \frac{|T_2|+1}{2}k + \epsilon + |U_1\cap U_2|(k+1)=\frac{|U_1|+1}{2}k + \epsilon + \frac{|U_2|+1}{2}k + \epsilon
+ |U_1\cap U_2|$, a contradiction. 

Combining the above two cases, we obtain (2). This contradiction justifies the claim.

\vskip 2mm

For each odd set $U \subseteq Z$, by the above claim, we obtain $|U|(k+1)=2|E(U)|+|F(U)|=|E(U)|+
|E^+(U)| \ge |E(U)|+\frac{|U|+1}{2}k +  \epsilon$. Thus $|E(U)|\le \frac{|U|-1}{2}(k+2)+1 -  \epsilon$.
Hence $\frac{2 |E(U)|}{|U|-1}\le k+3$ if $\epsilon=0$ and  $\frac{2 |E(U)|}{|U|-1}\le k+2$ if $\epsilon=1$.
By Conjecture \ref{GS}, the chromatic index of $G[Z]$ is at most $k+3 - \epsilon$. Since all vertices outside 
$Z$ have degree one, we further obtain $\chi'(G) \le k+3 - \epsilon$. So $E$ can be 
partitioned into $k+3- \epsilon$ matchings $M_1, M_2, \ldots, M_{k+3- \epsilon}$. 

Let us first consider the case when $\epsilon=0$. By the above claim, 

(7) each vertex $z\in Z$ is disjoint from precisely two of $M_1, M_2, \ldots, M_{k+3}$ (as $d(z)=k+1$).

Let $H$ be the subgraph of $G$ induced by edges in $M_k \uplus M_{k+1} \uplus M_{k+2} \uplus M_{k+3}$, where 
$\uplus$ is the multiset sum, and let $N$ be an orientation of $H$ such that $|d_N^+(v)-d_N^-(v)|\le 1$ for each 
vertex $v$. (It is well known that every multigraph admits such an orientation.) From (7) and this orientation we see that

(8) if a vertex $z\in Z$ is disjoint from precisely one of  $M_1, M_2, \ldots, M_{k-1}$, then $d_H(z)=3$
and $d_N^-(z) \ge 1$; if $z$ is disjoint from precisely two of  $M_1, M_2, \ldots, M_{k-1}$, then 
$d_H(z)=4$ and $d_N^-(z) = 2$. 

For each $i=1,2,...,k-1$, let $C_i$ be obtained from $M_i$ as follows: for each $z\in Z$, if $z$ not covered by $M_i$, 
add an edge from $N$ that is directed to $z$ and has not yet been used in $C_1 \uplus C_2 \uplus \ldots \uplus C_{i-1}$, 
where $C_0=\emptyset$. From this construction and (8) we deduce that $C_1, C_2,...,C_{k-1}$ are pairwise disjoint and 
each of them is a $Z$-cover in $G$. 

It remains to consider the case when $\epsilon=1$. Now

(9) each vertex $z\in Z$ is disjoint from precisely one of $M_1, M_2, \ldots, M_{k+2}$.

Let $H$ be the subgraph of $G$ induced by edges in $M_{k+1} \uplus M_{k+2}$, and let $N$ be an orientation of $H$ such 
that $|d_N^+(v)-d_N^-(v)|\le 1$ for each vertex $v$. From (9) and this orientation we see that

(10) if a vertex $z\in Z$ is disjoint from precisely one of  $M_1, M_2, \ldots, M_{k}$, then $d_H(z)=2$ and $d_N^-(z) = 1$. 

For each $i=1,2,...,k$, let $C_i$ be obtained from $M_i$ as follows: for each $z\in Z$, if $z$ not covered by $M_i$, 
add an edge from $N$ that is directed to $z$. From this construction and (10) we deduce that $C_1, C_2,...,C_{k}$ are 
pairwise disjoint and each of them is a $Z$-cover in $G$. \qed

\section{Implication}

The purpose of this section is to show that Conjecture \ref{CC2} can be deduced from Conjecture \ref{CC1}.

\vskip 2mm

{\bf Proof of Theorem \ref{bound}.} We may assume that

(1) $G$ is connected.

To see this, let $G_1, G_2, \ldots, G_k$ be all the components of $G$. For each $i=1,2, \ldots, k$, we aim to 
establish the inequality $\xi(G_i)>\delta(G)-\mu(G)$. If $\delta(G_i)-\mu(G_i) >\delta(G)-\mu(G)$, then the desired 
inequality holds, because $\xi(G_i)\ge \delta(G_i)-\mu(G_i)$. So we assume that $\delta(G_i)-\mu(G_i) \le 
\delta(G)-\mu(G)$. Since $\delta(G_i) \ge \delta(G)$ and $\mu(G_i) \le \mu(G)$, from this assumption we deduce 
that $\delta(G_i)=\delta(G)$ and $\mu(G_i)=\mu(G)$. Thus $G_i$ satisfies the hypothesis of Conjecture \ref{CC2}.  
Hence we may assume that $G$ is connected, otherwise we consider its components separately.  

By hypothesis, $\delta(G)$ cannot be expressed in the form $2p\mu(G)-q$, for any two integers $p$ and $q$ satisfying 
$q\ge 0$ and $p > \mu(G)+\lfloor (q-1)/2 \rfloor$; these two inequalities are equivalent to $0 \le q \le 2p-2 \mu(G)$. 
Setting $q=0, 1, \ldots, 2p-2\mu(G)$ respectively, we see that $\delta(G)$ does not belong to the set 
\[\Omega_p=\{2(p+1) \mu(G)-2p, 2(p+1) \mu(G)-2p+1, \ldots, 2p\mu(G)\},\]
where $p \ge \mu(G)$. Note that $2\mu(G)^2$ is the only member of  $\Omega_{\mu(G)}$ and that the gap between $\Omega_p$
and $\Omega_{p+1}$ consists of all integers $i$ with $2p\mu(G)+1 \le i \le 2(p+2) \mu(G)-(2p+3)$. So

(2) either $\delta(G) \le 2\mu(G)^2-1$ or $2p\mu(G)+1 \le \delta(G) \le 2(p+2) \mu(G)-(2p+3)$ for some $p \ge \mu(G)$. 

We may assume that $\delta(G)\ge 1$, for otherwise, $\delta(G)=2p\mu(G)-q$ for $p=q=0$, contradicting the hypothesis
of Conjecture \ref{CC2}. Thus $\mu(G)\ge 1$. 

To prove the theorem, it suffices to show that for any odd set $U$ of $G$, we have $\frac{2|E^+(U)|}{|U|+1}\ge 
\delta(G) -\mu(G)+1$, or equivalently,

(3) $2|E(U)|+2|F(U)| \ge (|U|+1)(\delta(G)-\mu(G)+1)$.

\vskip 2mm
Set $k=\mu(G)$ if $\delta(G) \le 2\mu(G)^2-1$ and set $k=p+1$ if $2p\mu(G)+1 \le \delta(G) \le 2(p+2) \mu(G)-(2p+3)$ 
for some $p \ge \mu(G)$. We consider two cases according to the size of $U$.

\vskip 2mm
{\bf Case 1.} $|U|\ge 2k+1$.

We divide the present case into two subcases.

{\bf Subcase 1.1.}  Either $U \subsetneq V$ or $U=V$ and $\delta(G)$ is odd. In this subcase,    

(4) $2|E(U)|+2|F(U)| \ge |U|\delta(G)+1$. 

\noindent Indeed,  if $U \subsetneq V$, then $|F(U)|\ge 1$ by (1). If $U=V$ and $\delta(G)$ is odd, then $G$ contains
at least one vertex of degree at least $\delta(G)+1$, because $|V|=|U|$ is odd and the total number of vertices with odd 
degree is even. Hence (4) is true.   

(5) $|U|\delta(G)+1 \ge (|U|+1)(\delta(G)-\mu(G)+1)$.

Note that (5) amounts to saying that $\delta(G)\le (|U|+1)(\mu(G)-1)+1$. If $\delta(G) \le 2\mu(G)^2-1$, then
$\delta(G) \le (2\mu(G)+2)(\mu(G)-1)+1=(2k+2) (\mu(G)-1) +1 \le (|U|+1)(\mu(G)-1)+1$. If $\delta(G) \le 2(p+2) 
\mu(G)-(2p+3)$, then $\delta(G) \le 2(k+1)\mu(G)-(2k+1) =(2k+2) (\mu(G)-1) +1 \le (|U|+1)(\mu(G)-1)+1$. So (5) is
established.

The desired statement (3) follows instantly from (4) and (5).

{\bf Subcase 1.2.}  $U=V$ and $\delta(G)$ is even. In this subcase, we have $\delta(G) \le 2\mu(G)^2-2$ if 
$\delta(G) \le 2\mu(G)^2-1$ and $\delta(G) \le 2(p+2) \mu(G)-(2p+4)$ if $\delta(G) \le 2(p+2) \mu(G)-(2p+3)$.
So $\delta(G) \le (2k+2)(\mu(G)-1)$ by the definition of $k$ and hence

(6) $\delta(G) \le (|U|+1)(\mu(G)-1)$.

From (6) we deduce that $|U|\delta(G) \ge (|U|+1)(\delta(G)-\mu(G)+1)$. Therefore (3) holds,
because $2|E(U)|+2|F(U)|\ge |U|\delta(G)$. 

\vskip 2mm
{\bf Case 2.} $|U|\le 2k-1$. (So $k \ge 2$ as $|U|\ge 3$.) 

By the Pigeonhole Principle, some vertex $v\in U$ is incident with at most $\frac{|F(U)|}{|U|}$ edges in
$F(U)$. Note that $v$ is incident with at most $(|U|-1) \mu(G)$ edges in $G[U]$, so $d(v)\le (|U|-1) \mu(G)+
\frac{|F(U)|}{|U|}$. Hence

(7) $\delta(G) \le (|U|-1) \mu(G)+ \frac{|F(U)|}{|U|}$.

We proceed by considering two subcases.

{\bf Subcase 2.1.}  $2p\mu(G)+1 \le \delta(G) \le 2(p+2) \mu(G)-(2p+3)$, where $p \ge \mu(G)$. 

From (7) and the hypothesis of the present subcase, we deduce that $2p\mu(G)+1 \le (|U|-1) \mu(G)+ 
\frac{|F(U)|}{|U|}$. Thus $|F(U)| \ge |U|(2p+1-|U|)\mu(G)+|U|$. So 

(8) $|U|\delta(G)+|F(U)| \ge |U|\delta (G) + |U|(2p+1-|U|)\mu(G)+|U|$. 

Let us show that 

(9) $|U|\delta (G) + |U|(2p+1-|U|)\mu(G)+|U| \ge (|U|+1)(\delta(G)-\mu(G)+1)$.

To justify this, note that (9) is equivalent to

(10) $\delta(G) \le \{|U|(2p+2-|U|)+1\}\mu(G)-1$. 

By the hypothesis of the present subcase, $\delta(G) \le 2(p+2) \mu(G)-(2p+3)$. To establish (10), we turn to
proving that $2(p+2) \mu(G)-(2p+3) \le \{|U|(2p+2-|U|)+1\}\mu(G)-1$, or equivalently

(11) $\{-|U|^2+2(p+1)|U|-(2p+3)\} \mu(G) \ge -(2p+2)$.

Let $f(x)=-x^2+2(p+1)x-(2p+3)$. Then $f(x)$ is a concave function on $\mathbb R$. So on any interval $[a,b]$,
$f(x)$ achieves the minimum at $a$ or $b$.  By the hypothesis of the present case, $|U|\le 2k-1=2p+1$, 
so $3\le |U|\le 2p+1$. By direct computation, we obtain $f(3)=4p-6\ge -2$ and $f(2p+1)=-2$. Thus $f(|U|)\ge -2$
for $3 \le |U| \le 2p+1$, which implies that the LHS of (11) $\ge -2 \mu(G) \ge -(2p+2) =$ RHS of (11), because 
$p \ge \mu(G)$. This proves (11) and hence (10) and (9).

Since  $2|E(U)|+2|F(U)| \ge |U|\delta(G)+|F(U)|$, the desired statement (3) follows instantly from (8) and (9).

{\bf Subcase 2.2}.  $\delta(G) \le 2\mu(G)^2-1$.  

We may assume that

(12) $\delta(G) \ge (|U|+1)(\mu(G)-1)+1$, for otherwise, $|U|\delta(G) \ge (|U|+1)(\delta(G)-\mu(G)+1)$. So
(3) holds, because  $2|E(U)|+2|F(U)| \ge |U|\delta(G)$. 

By (12) and the hypothesis of the present subcase, either $2t(\mu(G)-1)+1 \le \delta(G) \le 2(t+1)(\mu(G)-1)$ 
for some $t$ with $\frac{|U|+1}{2} \le t \le \mu(G)$ or $\delta(G)=2t(\mu(G)-1)+1 $ for $t=\mu(G)+1$. 
 
By (7), we have $2t(\mu(G)-1)+1 \le  (|U|-1) \mu(G)+ \frac{|F(U)|}{|U|}$. So 
$\frac{|F(U)|}{|U|} \ge (2t-|U|+1)\mu(G)-2t+1$, and hence

(13) $|U|\delta(G)+|F(U)| \ge |U|\{ \delta(G)+(2t-|U|+1)\mu(G)-2t+1 \}$. 

We propose to show that 

(14) $|U|\{ \delta(G)+(2t-|U|+1)\mu(G)-2t+1 \} \ge (|U|+1)(\delta(G)-\mu(G)+1)$.

To justify this, note that (14) is equivalent to  

(15) $\delta(G) \le \{ |U|(2t+2-|U|)+1\} \mu(G) -|U|2t-1$.

Suppose $\delta(G)=2\mu(G)^2-1$. Then $t=\mu(G)+1$. So (15) says that $2\mu(G)^2-1 \le \{ |U|(2\mu(G)+4-|U|)+1\} \mu(G) 
-|U|(2\mu(G)+2)-1$, or equivalently,  $\{ |U|(2\mu(G)+4-|U|)+1\} \mu(G) -|U|(2\mu(G)+2) \ge 2\mu(G)^2$.
Let $g(x)= \{ x(2\mu(G)+4-x)+1\} \mu(G)-x(2\mu(G)+2)$. Then $g(x)$ is a concave function on $\mathbb R$. So on any 
interval $[a,b]$, $g(x)$ achieves the minimum at $a$ or $b$. By direct computation, we obtain $g(3)= 6\mu(G)^2-
2\mu(G)-6$ and $g(2\mu(G)-1)= 6\mu(G)^2- 6\mu(G)+2$. It is easy to see that  $\min\{g(3), g(2\mu(G)-1)\} \ge 
2\mu(G)^2$, because $\mu(G)=k \ge 2$ (see the hypothesis of Case 2).  Hence $g(|U|) \ge 2\mu(G)^2$ for  
$3 \le |U| \le 2\mu(G)-1=2k-1$.  This proves (15) and hence (14) and (13).
 
So we assume that  $\delta(G) \le 2(t+1)(\mu(G)-1)$ for some $t$ with $\frac{|U|+1}{2} \le t \le \mu(G)$.  We prove
(15) by showing that $ 2(t+1)(\mu(G)-1) \le  \{ |U|(2t+2-|U|)+1\} \mu(G) -|U|2t-1$, or equivalently, 
$\{|U|(2t+2-|U|)-2t-1\}\mu(G)-|U|2t \ge -2t-1$. Let $h(x)= \{x(2t+2-x)-2t-1\}\mu(G)-2tx$. Then $h(x)$
is a concave function on $\mathbb R$. So on any interval $[a,b]$, $h(x)$ achieves the minimum at $a$ or $b$.
By direct computation, we obtain $h(3)=4(t-1)\mu(G)-6t$ and $h(2t-1)= 4(t-1)\mu(G)-2t(2t-1)$. It is easy 
to see that  $\min\{h(3), h(2t-1)\} \ge -2t-1$, because $\mu(G) \ge t \ge \frac{|U|+1}{2} \ge 2$. Hence 
$h(|U|) \ge -2t-1$ for  $3 \le |U| \le 2t-1$. This proves (15) and hence (14) and (13).

Since  $2|E(U)|+2|F(U)| \ge |U|\delta(G)+|F(U)|$, the desired statement (3) follows instantly from (13) and (14),
competing the proof of Theorem \ref{bound}. \qed\\

\end{document}